%% file: bstar.tex
\documentclass[a4paper,11pt]{amsart}
\input{commands}

\usepackage{hyperref}
\input xy
\xyoption{all}

\makeatletter
\renewcommand\paragraph{\@startsection
	{paragraph}{3}{0ex}
	{-\baselineskip}
	{-4ex}
	{\bf}}	
\makeatother

\newcommand{\formal}{[\![\lambda]\!]}
\newcommand{\laurent}{(\!(\lambda)\!)}
\newcommand{\lad}{\tfrac{i}{\lambda}\ad}
\newcommand{\fdef}[1]{\widehat{#1}}
\newcommand{\Weyl}{\mathcal{W}}
\newcommand{\bX}{\overline{X}}
\newcommand{\stimes}{\!\times\!}

\newcommand{\op}{^{\mathrm{op}}}
\newcommand{\A}{\mathcal{A}}
\newcommand{\C}{\mathcal{C}} 
\newcommand{\I}{\mathcal{I}}
\newcommand{\bI}{\boldsymbol{\mathcal{I}}}

\begin{document}

\title{Bohr-Sommerfeld Star Products}

\author{Michael Carl}

\address{
Fakult\"at f\"ur Mathematik und Physik,
Albert-Ludwigs Universit\"at Freiburg,
Hermann Herder Str. 3,
79104 Freiburg,
Germany} 

\email{michael.carl@physik.uni-freiburg.de} 

\begin{abstract}
We relate the Bohr-Sommerfeld conditions to formal deformation quantization of symplectic manifolds by 
classifying star products adapted to some Lagrangian submanifold $L$, 
 i.e. products preserving the classical vanishing ideal $\bI_L$ of  $L$ up to 
$\bI_L$-preserving equivalences.
\end{abstract}

\maketitle

\section{Introduction and Motivation}

 \paragraph{Reminder on Bohr-Sommerfeld conditions}
 Let $L\stackrel{i_L}{\inj} T^*Q$ be a Lagrangian submanifold of some cotangent bundle $T^*Q\stackrel{\pi}{\to}Q$ with respect to the standard symplectic structure $\omega:=-d\theta$ given by the canonical form $\theta :=T^*\pi$,
 and $\mu$ its  Maslov class. 
 Then $L$ has to satisfy
 the prequantization condition
 \begin{equation}\label{BS}
   \tfrac{1}{2\pi \lambda}i^*_L\theta  - \tfrac{\pi}{2} \mu 
\in  H^1_{dR}(L, \mathbb{Z}). 
 \end{equation}
 in order to be the microsupport of some $\lambda$-oscillatory distribution on $Q$
 (see Appendix for references), where $  H_{dR}(., \mathbb{Z})$ denotes the integral de Rham classes.
 
 In particular, if $L_E:=H^{-1}(E)$ is a Liouville torus of some semiclassical integrable system
with classical moment map $H: T^*Q\to \real^n$ and  vanishing subprincipal form\footnote{Let $s_i$ be the second order part of the Weyl symbols of the $\hat{H_i}$, then by definition $\kappa .\ham H_i = -s_i$.} $\kappa$,
then the condition \eqref{BS} coincides up to higher orders $O(\lambda^1)$ with 
the Bohr-Sommerfeld conditions
 \begin{equation}\label{BS2}
    \tfrac{1}{2\pi \lambda}i^*_L\theta  - \tfrac{\pi}{2} \mu  + \kappa +O(\lambda) \in  H^1_{dR}(L, \mathbb{Z})
    \end{equation}
 for the existence of a joint asymptotic eigenvector  of  the system. 

Recall here that a semiclassical integrable system
 is a maximal set of commuting $\lambda$-pseudodifferential operators $\hat{H}_1,...,\hat{H}_n$
 whose principal symbols $H=(H_1,...,H_n)$ are independent almost everywhere. 
 Then the compact fibers $L_E$ 
 of $H$ restricted to regular values $B$ are Lagrangian submanifolds with a transitive locally free  $\real^n$ action $\phi: \real^n\ni t \mapsto\exp \ham \left<t,H\right>$
 generated by the Hamiltonian vector fields $\ham H_i:=\omega^{-1}dH_i$, 
hence they are tori 
 $L_E\cong \real^n/\ker \phi|_{L_E}  \cong \real^n/\mathbb{Z}^n$.
This action is linearized in action-angle coordinates;
 similarly, the semiclassical system is 
 microlocally unitary equivalent to
 the linearized system, such that the microlocal solutions (oscillatory constants) form a local system
 whose triviality is the Bohr-Sommerfeld condition \eqref{BS2}, cf.
\cite{gnoc:1998a}.

\paragraph{Problems with interpretations in deformation quantization} 
The  Bohr-Sommerfeld condition \eqref{BS2} makes precise the original notion of quantization as ad hoc-discretization of classical spectra. However, 
it is undefined in formal deformation quantization, a notion isolating 
the transition from 
commutative to non-commutative algebras  
underlying any quantization concept (see Appendix \eqref{eq:*prod} for references). Here the
deformation parameter $\lambda$ is formal, one thus has
consider convergent deformations (of a subalgebra) over some base including $[0,\hbar]$
to recover the meaning of $\eqref{BS}$.
 However, we can try to extract the prequantization class in \eqref{BS} from deformation
 quantization as formal class. This has been done in \cite{nest:2004}
 by formalizing the symbol calculus of oscillatory distributions. Here, we
 proceed differently:

\paragraph{Main results and outline}
Motivated by oscillatory symbol calculus as well, we first consider $\star$-\-re\-pre\-sen\-ta\-tions on
line bundles over $L\subset X$. Then we look for star products inducing a canonical
representation on $L$, namely, 
we establish a bijection between the classes of
deformation quantizations of $X$ preserving the classical vanishing ideal $\bI_L\subset
C^\infty(X)\formal$ of $L$ up to $\bI_L$-preserving  equivalences
and those of formal deformations of the symplectic form viewed as relative class.
This is done in section \ref{sec:adapted} by parametrizing  adapted Fedosov star products. 
Then in section \ref{sec:relations} we consider the induced intertwiners on the quotients
as formal analogues of $L$-deformations in order to
explain the coincidence of \eqref{BS} with adapted classes in lowest order.
Finally, we sketch relations to the Maslov index and the symbol calculus of
oscillatory distributions.
The Appendix provides the motivating background
and sets up some (standard) notations.

\section{Representations on line bundles}\label{sec:rep}
Let $\star$ be a star product on a symplectic manifold $X$ and $E$
some vector bundle over a Lagrangian submanifold $L \subset X$.
A $\star$ representation on  $E$ means a $\star$-module structure on
$\Gamma(E)\formal$
such that $\star$ acts by $\comp\formal$-linear differential operators. 

Finally, define a deformed flat line bundle on $L$ by the requirement that its sheaf of local sections is locally isomorphic to the 
constant sheaf $\exp (\comp \formal)_L$ of the units 
$\exp (\comp \formal) \subset \comp \formal$. The isomorphism classes of such bundles are then given by
$\check{H^1}(X,\exp (\comp \formal)_L)$ like in the undeformed case of flat line bundles.
This group now acts naturally on $\star$ representations on $E$ 
thanks to their $\comp\formal$-linearity. The action turns out to be free and transitive:
\begin{lemma}\label{lemma:rep}
$\star$ is representable on some complex line bundle $E$ over a Lagrangian submanifold $L$ if and only if
the image $c_1^{\real}(E)$ of its Chern class $c_1(E)$ under $^{\real}: H^2(L,\mathbb{Z})\to H^2_{dR}(L)$ coincides with the restriction of  the equivalence class $[\star]$ of $\star$:
\begin{equation}\label{intc}
c_1^{\real}(E)=i_L^* [\star].
\end{equation} 
The space $\mathcal{M}_E$ of isomorphism classes of $\star$ representations on $E$ identifies with the group
of deformed flat line bundles $\check{H}^1(L, \exp({\comp \formal})_L)$ via its natural action on $\mathcal{M}_E$.
\end{lemma}
\begin{proof}
This is a direct consequence of Bordemann's classification:
Consider the restriction of $\star$ to some
tubular neighborhood of $L$ which we may identify with a 
neighborhood $W$ of the zero section in $(T^*L, -d\theta)$  
by some Weinstein isomorphism (cf. \cite[Th. 4.19]{bates.weinstein:1995a}).
Then by \cite[Th. 3.3]{bordemann:2004} the product $\star|W$
is equivalent to a standard ordered product $\star_B$, $[B]=i_L^*[\star]$,
whose
representation $\bullet$
on some complex line bundle $E$ on $L$ 
must locally on some contractible set $U_i$ look like 
\begin{equation}\label{connection}
(\pi^*\psi \bullet \phi) = \psi \phi , \qquad   (\theta.\hat{X}) \bullet \phi |U_i = -2\lambda (X - \tfrac{A_i.X}{2\lambda} ) \phi|U_i
\end{equation}  
for any
$\phi \in C^\infty(L)$ and any  
vector field $X\in \Gamma (TL)$ with canonical lift $\hat{X}\in \Gamma (T_0TL)$. Here the $A_i$ are determined by  $dA_i= B|U_i$ up to some coboundary $dS_i$, which determines $\bullet |U_i$ 
up to some local intertwiner (gauge equivalence) $\phi|U_i \mapsto e^{S_i}\phi|U_i$. 
Thus \eqref{intc} must hold,
which determines the representation up to 
isomorphism classes of formal flat connections 
$$ 
\lambda H^1_{dR} (L)/2\pi i H^1_{dR} (L;\mathbb{Z}) + \lambda^2    H^1_{dR} (L)\formal $$
on the chosen torsion bundle in $\ker ^{\real}$. 
\end{proof}

\begin{remark}
\newcommand{\F}{\mathcal{F}}
Consider two line bundles $E,E'$ over $L$. If $E$ is a $\star$ representation 
and $c_1(E')-c_1(E)\in \im(i^*_L: H^2(X,\mathbb{Z})\to H^2(L,\mathbb{Z}))$
then one may obtain
a $\star'$ representation on $E'$ via Rieffel induction:

Namely, it was shown in \cite{bursztyn.waldmann:2002a} that 
$\Pic(X)\ltimes \Aut(X,\omega)$ acts transitively on the Morita equivalent
equivalence
classes of $\star$
such that 
$[(\mathcal{L},\phi) \cdot \star] = [\star] + c_1^\real(\mathcal{L})$.
Here the symplectomorphism $\phi\in \Aut(X,\omega)$ acts by pull back while the action of 
the line bundle
$\mathcal{L}$ on $X$ arises by 
deforming its transition 1-cocycles to $\star$-cocycles defining an
equivalence $(\star',\star)$-bimodule $\boldsymbol{\mathcal{L}}$, cf. 
\cite[sec. 4.2]{bursztyn.waldmann:2002a}.
Hence $\boldsymbol{\mathcal{L}} \otimes_{\star} E$ indeed defines a $\star'$ representation on $i^*_L{\mathcal{L}} \otimes  E$.

\end{remark}

\section{Adapted star algebras}\label{sec:adapted}

Recall that a star product $\star$ on $X$ is called adapted to some Lagrangian submanifold $L\subset X$ if the classical vanishing ideal $\I_L:=\set{f\in \A}{ f|L=0}$ of $L$ in $\A:=C^\infty(X)$ remains a $\star$-left ideal
$\bI_L=\I_L\formal$ 
and thus induces a $\star$ representation
$\A\formal /\bI_L $ on $L$. Such products are formal series of one chains of the 
subcomplex
$K_{\I_L}:= \set{C\in \C^\bullet(\A;\A)}{ C( \A^{\otimes \bullet} \otimes \I_L) \subset \I_L} $
of the differential Hochschild complex\footnote{
Recall that the differential Hochschild complex $\C^\bullet(\A,M)$ of $\A$ with values in some 
$\A\otimes \A^{\op}$-representation $M$ consists
of $M$-valued differential $\comp$-multilinear operators $\C^{k-1}(\A,M):=\Hom_{\comp}^{\mathrm{diff}}(\A^{\otimes k}, M)$ with coboundary 
$bC(f_0,...,f_n) \!=\!
 f_0 C(f_1,...,f_{n}) \!+ \! \sum_i (-1)^{i+1} C(f_0,...,f_if_{i+1},...f_{n})+
(-1)^{n-1}C(f_0,...f_{n-1})f_{n}$.
} 
\newcounter{foot}
\setcounter{foot}{\value{footnote}}
of $\A$ fitting into an exact sequence
$$ K_{\I_L}\inj  \C \sur \C(\A;\C^1(\I_L;\A/\I_L))[-1] $$
by \cite[Prop. 2.2]{bordemann.ginot:2005}.   
Its corresponding long exact cohomology sequence decouples into short exact sequences isomorphic to the one
defining relative de Rham forms
\begin{equation}\label{reldR}
\Omega(X,L) \inj \Omega(X) \stackrel{i_L^*}{\sur} \Omega(L).
\end{equation}
This was shown locally in \cite[Th. 2.4]{bordemann.ginot:2005} for the $\omega$-corresponding multivector fields 
via Koszul resolutions,
from which the global case can be deduced by the degeneration of the local-to-global spectral sequence at $E^{pq}_2= H^p(\Omega^q(X,L)_X)= \Omega^q(X,L)\delta_{p0}$.
However, 
we will not use this ``adapted HKR theorem", although together with \eqref{Lichnerowicz1} it implies directly 
the following Lemma: Denote by $\delta$ the connecting homomorphism of the long exact sequence 
 \begin{equation}\label{eq:longexact} 
 ... H^\bullet_{dR}(X)\to H^\bullet_{dR}(L) \stackrel{\delta}{\to}   H^{\bullet +1 } _{dR}(X,L)\to H_{dR}^{\bullet +1}(X)... 
 \end{equation} 
associated to
the relative de Rham sequence \eqref{reldR}. Then we have:

\newcommand{\can}{\mathrm{can}} 
\begin{lemma}\label{AdaptedEquivalence}
Let $S$ be an equivalence between two star products $\star, \star'$ on $X$ both adapted to $L$. 
If $\delta H^1_{dR}(L) =\{0\}\subset   H^2_{dR}(X,L)$,
then $S$ is adapted, i.e. preserves $\bI_L$ 
and hence provides an equivalence\footnote{
\renewcommand{\H}{\mathcal{H}}
Following \cite[Prop 3.4]{bordemann:2004}
we call two represented algebras $(\A,\H),(\A',\H')$ equivalent if there is an
isomorphism $S: \A \to \A'$ and an
 ``intertwiner"
$T:\H\bij \H'$ such that $ T( f\cdot \psi) =S(f) \cdot   T\psi $.
}
of represented algebras 
$(\star,\star/\bI_L)\sim (\star',\star'/\bI_L)$.
\end{lemma}
\begin{proof}
Suppose 
$\star, \star'$ are already identical up to $O(\lambda^{k})$
thanks to an equivalence adapted up to $O(\lambda^k)$. Then by \eqref{Lichnerowicz1} and \eqref{HKR}
one has 
Lichnerowicz's equation
\begin{equation}\label{Lichnerowicz}
 \star' - \star  = bT + d\alpha(\ham .,\ham.) \mod O(\lambda^{k+1}) 
\end{equation}
 for some 
$\alpha \in \Omega^1(X)$, 
where $b$ is the Hochschild coboundary of the undeformed product '$\cdot$'. 
As decomposition into symmetric and antisymmetric parts, both summands have to be adapted by induction hypothesis, so $d i^*_L\alpha =0$. 
Now by the assumption 
there exists  some relative primitive of $d\alpha$
we may suppose to be $\alpha$ itself. 

Then the equivalence 
$S_\alpha:=1+ \lambda^{k-1} \alpha \ham $ is adapted
and turns the difference into a coboundary
$$S_\alpha(\star') - \star   = \lambda^k b(T - [T,\alpha \ham]) =: \lambda^k bT' \mod O(\lambda^{k+1})$$
which is adapted if and only if $T'$ is, since $bT'(\I_L,\I_L) \subset  T'(\I_L^2) + \I_L$.
Thus $S':= (1 + \lambda^k T')\circ S_\alpha$  is an adapted
equivalence modulo $O(\lambda^{k+1})$. 
\end{proof}

 \paragraph{Reminder on Fedosov's construction}
\newcommand{\fweyl}{\fdef{\real^{2n}}} 
\newcommand{\derweyl}{\mathfrak{g}}
There is a construction of natural\footnote{Following
\cite{gutt.rawnsley:2003a}, we call a star product $\star$ natural if  $\star_k$ is of differentiation order $\leq k$ in each argument.}
  Weyl type star products due to Fedosov \cite{fedosov:1996a} which 
has a natural interpretation in formal geometry 
context as observed by \cite{nest.tsygan:1995b}:

Let $\fweyl:=  (\comp [\![ \xi^1,...,\xi^{2n};\lambda]\!],*_W)$
 denote the formal Weyl algebra 
 \begin{equation}\label{weyl} f*_Wg:= \mu_0(\exp({\tfrac{\lambda}{2i} \omega ^{ij} \partial_i \otimes \partial_j}) f\otimes g ),
\end{equation}
where $\mu_0$ denotes the standard multiplication. $*_W$ 
is $\mathbb{Z}$-graded by
$\deg \lambda =2$, $\deg \xi = 1$
and invariant under the 
linear symplectic group $ Sp(n,\comp) \cong  \ad \fweyl_2 $  generated by quadratic forms in $\fweyl_2$, where as usual the subscript denotes the degree.

One now considers infinitesimal patching of local algebras on $X$, i.e. the bundle $P$ of
isomorphisms $jet_x(X)\formal  \bij \fweyl$, $x\in X$, deforming the real symplectic frame bundle 
$Sp(X)= \bigcup_x (T_xX,\omega_x)\bij  (\real^{2n}, \omega_0)$. Then the natural 
isomorphism $\theta: T_xP\bij \derweyl:=\set{\lad f}{f\in \fweyl,  \Im(f)=0\mod O(\lambda)} $ is a flat $\derweyl$-valued connection, i.e. a $Sp$-equivariant
1-form such that its composition with the $\mathfrak{sp}$-action is the inclusion $\mathfrak{sp} \to \derweyl$ and
the curvature $d\theta +\half [\theta,\theta]$ vanishes.
Now any reduction of the structure group given by a section $r$ of $P\to Sp(X)$
induces a flat connection $\nabla_F = d + \lad r^*\theta$  on the associated Fedosov bundle
   $\Weyl:= Sp(X)\times_{Sp(n)} \fweyl $. 
           By \cite{nest.tsygan:1995b}, its constant sections $\ker \nabla_F$ are isomorphic to a star product algebra on $X$ via $r$ whose
characteristic class is represented by the pullback by $r$ of the curvature
   of the lift of $\theta$ to a connection with values in the central extension $\fweyl$ of $\derweyl$. 
   
   More explicitly, any Fedosov connection $\nabla_F$ has to start with an equivariant degree $-1$ square zero
   differential fixed as
   $$-\delta:=-\lad r^*\theta_1:= -\tfrac{\partial}{\partial \xi^i}\otimes dx^i.$$ Further the $\lambda$-independent degree zero part $\base{x^l} + \lad \Gamma ^{l}_{jk} \xi^j\xi^k$
represents a torsion free symplectic connection $\nabla$.
   Thus the simplest Fedosov connection is of the form
   $ \nabla_F= -\delta + \nabla +\lad_{*_W}{ \gamma} $
 for some $\gamma \in \Omega^1(X; \Weyl_{\geq 3})$.

However, to get natural products of different order type, one needs
equivalent fiberwise degree 0 products $* = \mu_0 \exp(\frac{\lambda}{2i}\mu_{ij}\base{\xi_i}\otimes \base{\xi_j})  = e^{\lambda S}(*_W)$, where 
$S$ is a fiberwise degree $-2$ differential operator with fiberwise constant coefficients, which in turn
requires Fedosov derivations of the generalized form (cf. \cite{neumaier:2001a})
\begin{equation}\label{eq:FedDer}
 \nabla_F = -\delta + D + 
\lad{\gamma},
\end{equation}
where $D$ is a degree 0 $*$-derivation with $D|1\otimes \Omega(X)=d$, 
$ D^2 = - \lad R$ and $[D,\delta] = \lad T$ for totally covariant curvature and torsion
tensors $R$ and $T$. Namely, we may and will take the torsion free derivation
\begin{equation}\label{eq:D}
D= \nabla + \tfrac{i\lambda}{2} [ \nabla, S].
\end{equation}
Note that at this point it might become more natural to work with
the deformation of the symmetric algebra $ST^*X$ 
induced by $(\Weyl, *)$ and the canonical isomorphism
$\Weyl /\lambda \Weyl \cong ST^*X$, as $(\Weyl,*)$ is no longer given as associated bundle.

 Now,  for any formal closed two form $\Omega\in Z^2_{dR}(X)\formal$ and any $s\in \Weyl_{\geq 4}$
 with trivial central part $\sigma(s) = 0$ 
 the equations 
 \begin{equation}\label{gamma} 
 \delta \gamma = D \gamma - \tfrac{1}{\lambda} \gamma * \gamma + R + \Omega, \quad \delta^{-1}\gamma = s
 \end{equation}
 determine $\gamma$ such that $\nabla_F$
 will be a Fedosov connection whose constant sections $(\ker \nabla_F,*)$
are naturally isomorphic to a natural star product $\star_F$ of class 
$[\star_F]=[\lambda^{-1}\omega +\Omega]$. More precisely, $\nabla_F$
extends to an acyclic superderivation on $\Omega(X,\Weyl)$ with contracting 
homotopy 
$$\nabla_F^{-1}\alpha := -\delta^{-1} \frac{1}{1-[\delta^{-1}, D+\lad(\gamma) ]}$$ 
(in the sense of the geometric series), where $\delta^{-1}$ is defined as homotopy
on the center
$(\delta^{-1}\delta + \delta \delta^{-1}) \alpha  =  \alpha -\sigma(\alpha)$.
Then the isomorphism 
$\star_F \cong  \ker \nabla_F \cap \Weyl$ 
is given by the restriction of
$\sigma$ with inverse $\tau(f):= f-\nabla_F^{-1}df$. 
 
Using \cite[1.3.25,27]{neumaier:2001a}, one checks that the parametrization by $\nabla, \Omega, *, s$ is redundancy free for fixed $*$. Moreover, it is conjecture that any natural star product arises as generalized Fedosov star product.

\begin{proposition}\label{adaptFed}
A generalized Fedosov star product $\star_F$ on $T^*L$ 
given by (\ref{eq:FedDer},\ref{eq:D},\ref{gamma})
is adapted to $L$
if and only if its construction data $\nabla, \Omega, *, s$ are adapted, i.e:
\begin{enum}
\item $\nabla$ restricts to a connection on $L$. By absence of torsion this is equivalent to $L$ being totally geodesic.
\item  $L$ is Lagrangian for the ``deformed symplectic structure" $\omega + \Omega $
\item $s\in I_{TL}$, where $I_{TL}= (TL^0) \formal $ is the ($\comp \formal$-extended) fiberwise vanishing ideal of $TL$ generated by the annulator $TL^0\subset ST^*X$ of $TL$. 

\item The fiberwise product $*$ is adapted to $I_{TL}$, i.e., in symplectic fiber coordinates
$q^i, p_j$  over $L$ such that $I_{TL}=(p_1,...,p_n)$ we have $S=\frac{\lambda}{2i}\base{p_j}\base{q^j}$ and $* = \mu_0 \exp(\frac{\lambda}{2i}\base{p_j}\otimes \base{q^j}).$

\end{enum}
\end{proposition}
\begin{proof}
Let $I_F^{s,p} \subset \Omega^p(X,\Weyl_s)$ be the subspace of adapted forms whose restriction to
$\wedge^p TL$ has per definition values in $I_{TL}$. 
Now, if all construction data are adapted, then $\gamma \in I_F$,
since the same holds for all summands in \eqref{gamma}. Here
 the only non obvious term is $[\nabla, S]$, where the claim holds if 
 $\nabla$ is the homogenous adapted connection $\nabla^0$ used in \cite{bordemann.neumaier.waldmann:1998a}, then it follows in general by
 $\nabla - \nabla^0 \in \lad I_F^{2,1}$.
 Thus $\nabla_F$ preserves $I_F$ such that $\tau \bI_L= \ker \nabla_F \cap I_{TL}$,
 hence $\star_F$ is adapted. \---
 Vice versa:

{\em i.} Let $\sigma'$ denote the projection 
$\Omega(X,\Weyl) \to \Omega(X,\comp)\formal$.
We have $\sigma' \delta \tau f=  \delta \delta^{-1} df = df$.
Now $df|TL=0$ for any $f\in \bI_L$, 
so adaptivity
implies
\begin{equation}\label{eq:ideal}
\sigma'  \delta(X * \tau I)|TL=0 \quad \forall X\in \ker \nabla_F, I\in \bI_L
\end{equation}
with 
$\tau f = f + Df + D^2f$ modulo $\deg \geq 3$, where $D=[\delta^{-1},\nabla]$ 
here and in the following denotes
 the induced symmetric covariant connection on $\Weyl$.
Thus \eqref{eq:ideal}  implies 
\begin{equation}\label{eq:X1}
0= (X *_1 \delta D^2I)|TL 
\end{equation}
for some  $X$ of total degree 1.
By iv.
$X *$ only differentiates along $L$, thus for any
vector fields $X,Y$ 
tangential to $L$ one has 
$ 0 = D^2I(X,Y)|L =  dI(\nabla_XY + \nabla_YX)|L = dI( 2\nabla_X Y)|L$ by absence of torsion and $dI.[X,Y]|L=0$,
thus $\nabla_XY$ must be tangential to $L$. 
\\ {\em iv.} Since $\tau^k = D^k $ modulo lower order differentiation along $Z$, we have
 \begin{equation}
 f (\star_{F})_k I = \mu_{IJ} \partial^I f \partial^J I \mbox{ modulo lower order differentiation}   
 \end{equation}
 for some non degenerate tensor $\mu \in \Gamma (\Weyl_k \otimes \Weyl_k)$. In particular, 
 adaptivity and \eqref{eq:defLie} imply $\mu_{ij} = \base{p_j}\otimes \base{q^j}$ for $k=1$.  
\\ {\em ii. and iii.} Allow $s\in \Weyl_{\geq 3}$ to make $\Omega$ locally redundant as in 
\cite[1.3.27]{neumaier:2001a}, and denote by
$\tau(s)$ the dependence of $\tau$ on $s$. Then by Lemma \ref{AdaptedEquivalence} and adaptivity of
$\star_{F}$ for adapted construction data including $s=0$  
we must have $ \sigma( \tau(s)^{2k}-\tau(0)^{2k})\bI_L\subset \bI_L$.
Now  $s^k$ first occurs deg-inductively  in $\tau^{2k-2}$ as  
 $[\delta s^{k}f , \tau^{k-1} I]_{*_{k-1}}$, hence 
 $  \tau^{k-1}  I *_{k-1} \delta s^{k}f | TL = 0 $ for all $I\in \bI_L$, thus $s$ is adapted by standard order of $*$.
\end{proof} 
 
\begin{remark} \label{rem:Omega}
 Note that $\Omega_{k}$ enters explicitly as
 $\star_{\nabla, \Omega + \lambda^k\Omega_k, \circ, s} = \star_{\nabla,\Omega,\circ, s} +  
 \lambda^{k+1}\Omega_k(\ham ., \ham .) \mod O(\lambda^{k+2}) $.
\end{remark}

 Now identify some zero section environment of $T^*L$ with some neighborhood of 
 $L\subset X$ via some Weinstein isomorphism. Then note that there are no obstructions
 for the extension of the data to all of $X$: For the fiberwise data this is clear
 from partition of unity, for the connection this follows
 from the description of the space of symplectic connections as
 sections of the fiber bundle $J^1 Sp(X)/Sp$ having contractible fibers
 $F:=J^1_0(\real^{2n},Sp)_1 = \real^{2n} \times \mathfrak{sp} $ (see \cite{michor:1993}),
 hence the obstruction classes $H^i(X,L;\pi_{i-1}(F))$ vanish. 
 (Note that this argument ignores torsion, which is possible due to \cite[prop 1.3.31]{neumaier:2001a}).

 Now by \eqref{intc} 
the characteristic class $[\star]$ of an adapted star product must have a relative representative 
 $\eta \in Z^2_{dR}(X,L)\laurent$. For any relative cocycle $\eta$ 
 the above construction yields an adapted Fedosov star product
 $\star_F$ with $\eta= \lambda^{-1}\omega +  \Omega$, and two adapted Fedosov products differing only in their 2-forms by
 $\lambda^{k} d\alpha \mod O(\lambda^{k})$ are equivalent through 
$S_\alpha:= 1 + \lambda^{k-1}\alpha.\ham \mod O(\lambda^{k})$ by Remark \ref{rem:Omega} and \eqref{Lichnerowicz},
which covers the cohomology in \eqref{Lichnerowicz}.
 Since $S_\alpha$ is
 adapted if and only if $i_L^*\alpha =0$, i.e. 
 $[d\alpha] = 0\in H^2_{dR}(X,L)$,
 we obtain: 
 \begin{theorem}\label{th:classification} 
 The adapted equivalence classes $[.]_L$ of  $L$-adapted deformation quantizations are in bijection to relative formal $\tfrac{\omega}{\lambda}$-deformations
 $$\lambda^{-1}[\omega] + H^2_{dR}(X,L)\formal$$
 where the image and kernel of $H^2_{dR}(X,L)\formal$ in the long exact sequence \eqref{eq:longexact} correspond to
 the absolute class and the adapted classes therein, respectively.

 \qed
 \end{theorem}
 
 \begin{remarks}
 \item  
In the preprint \cite{bordemann.ginot:2005} the generalized problem of
deforming the adapted HKR quasiisomorphism to an adapted $G_\infty$ resp. $L_\infty$ morphism
has been attacked locally, see as well \cite{cattaneo.felder:2005a}.
 
\item  As shown by Gromov \cite[Th. 7.34]{mcduff.salomon:1995}, for an open manifold any de Rham cohomology class is representable by a symplectic form unique up to isotopy. This of course doesn't hold in the relative case, for instance  the trivial class in $H^2_{dR}(\real^2, S^1)$ has no symplectic representative $\omega$ by
 $\int_{\mathrm{int}S^1} \omega\not = 0$. 
In fact,
 the same holds for any compact Lagrangian submanifold $L\subset \real^{2n}$ by the above isotopy theorem and Gromov's nonexistence theorem of exact $L\subset \real^{2n}$ \cite[Th. 13.5]{mcduff.salomon:1995}.

\item By the same method, we can construct star products adapted to transversal  intersections of Lagrangian submanifolds or all the fibers of  a regular Lagrangian fibration.

 \end{remarks}

\section{Relation to Bohr-Sommerfeld conditions} \label{sec:relations}

We now want to ``deduce" the prequantum Bohr-Sommerfeld conditions from 
the picture
\begin{equation}\label{diagramm}
 \xymatrix@R=3.2ex{
\displaystyle\frac{\tfrac{1}{\lambda}i^*_L\theta + H^1_{dR}(L)\formal}{ \im \tfrac{1}{\lambda} H^1_{dR}(X)\formal} \; \ar@{^{(}->}[r]^{\delta} \ar[d]\ar@{<->}[d]
&
\tfrac{\omega}{\lambda} + H^2_{dR}(X,L)\formal \ar[r]\ar@{<->}[d]
&
\tfrac{\omega}{\lambda} + H^{2 } _{dR}(X)\formal\ar@{<->}[d]
\\
\txt{ adapted classes \\inside an absolute class}   
&
\txt{adapted classes }
&
\txt{ absolute classes}
  }
 \end{equation}
following from theorem \ref{th:classification}.

\begin{lemma}\label{equiIntertwiner}
Consider the group of equivalences $S_\alpha=e^{\lambda \alpha.\ham}$ between adapted star products modulo adapted equivalences, which is ho\-mo\-mor\-phi\-cal\-ly pa\-ra\-me\-trized by 
$i^*_L[\alpha] \in H^1_{dR}(L)\formal$. Then the  action of $S_\alpha$
on equivalence classes of pairs {\rm \{(adapted star product, canonical representation)\} }induces the action of the flat deformed line bundle $E_\alpha$ with 
holonomy $\pi_1(L) \ni \gamma \mapsto e^{i\lambda \int_\gamma \alpha}$ on representation classes
$\mathcal{M}_0$ (Lemma \ref{lemma:rep}).
\end{lemma}
Indeed, in a Weinstein model $T^*L$ near $L$ with Darboux coordinates $I,\varphi$  we can assume 
$\alpha = \alpha_i d\varphi_i$ near $L$ for $\alpha_i \in \comp\formal$, then $S_\alpha$ acts on contractible patches like the inner automorphisms  $S_\alpha = e^{\lambda \alpha_i \{ \varphi_i,\cdot\} }
=e^{\lambda \alpha_i \ad_\star(\varphi_i)} $ on the standard orderd product $\star$ of $T^*L$. As inner automorphisms induce self-intertwiners of $\star/\bI_L$, globally we get the desired action of $E_\alpha$ on 
$\mathcal{M}_0$. \qed
  
\smallskip

Hence in case of convergence, {\em  integral} adapted equivalences $S_\alpha, \lambda i^*_L\alpha\in H^1_{dR}(L,\mathbb{Z})$ should provide intertwiners 
between formally in-equivalent representations which correspond to relative quantization conditions as follows:

\paragraph{Relative conditions}
\newcommand{\reg}{\mathrm{reg}}

Restrict $H$ to a regular torus fibration
$H_{\reg}=H|X_{\reg}$ over $B$ 
near $L$.
Then there is a natural identification
$\mathbf{R}^1H_{\reg} \real = T_{B}$ obtained as dual of the action differential
$\gamma_i \mapsto d I_i:= d \int_{\gamma_i} \theta $.

In particular, we can {\em canonically} identify small $\alpha_i[d\varphi_i] \in H^1_{dR}(L)$ with the
translated torus $L_\alpha:=H^{-1}_{\reg}(\alpha_1,...,\alpha_n)$ (in action coordinates with origin $0=H(L)$) which equals the image $ \im \alpha_i d\varphi  \subset T^*L$ in the Weinstein model near $L$ determined by the action-angle coordinates.

Now $S_\alpha$ maps the vanishing ideal $\I_L$ of $L$ to that of the translated Lagrangian torus 
$L_\alpha$: In fact, infinitesimally, this map corresponds to the isomorphism
\begin{eqnarray*}
 \Hom_{\A}(\I_L/I^2_L, \A/\I_L)& \cong  \Gamma(T_LX/TL) &\cong \Omega^1(L): \\  \tdiffo{\lambda} S_\alpha = \alpha.\ham & \mapsto \quad [\omega^{-1}\alpha|L] &\mapsto i^*_L\alpha 
 \end{eqnarray*}
where $\A:= C^\infty(X)$, $\I_L:= \ker i_L^* \subset \A$ are the classical $\A$-modules.

Hence for integral $i^*_L\lambda \alpha$, $S_\alpha$ intertwines the canonical $\star$ representations ($D$-modules)
$\star/\I_L, \star/\I_{L_\alpha}$ on $L$ and $L_\alpha$.
But such intertwiners correspond to joint $(\alpha_1,...,\alpha_n)$-eigenspaces of $H_{\reg}$ by
the standard $D$-module identity (cf. \cite[Ch. 0]{kashiwara:2000})
\newcommand{\homD}{\hom_{D}}
\renewcommand{\L}{\tfrac{\A}{ \I_L}} 
 $$
 \ker \left( (H_{\reg,i}   - \alpha_i): \L \to \left(\L\right)^n \right) =       \homD\left(\tfrac{D}{(H_{\reg,i} -\alpha_i)}, \L\right).  $$

\begin{remark}\label{relAbs}
Note that our reference Lagrangian $L\subset X_{\reg}$ itself is always quantizable 
and in order to speak of non trivial relative classes on $X_{\reg}$, one has to identify the class 
$\alpha_i[d\varphi_i]\in H^1_{dR}(L)$ with 
$\alpha_i[d\varphi_i] \oplus 0 \in H^1_{dR}(-L\cup L_\alpha)$ giving a meaningful class in
$H^2_{dR}(X_{\reg},-L\cup L_\alpha)$.Then by \eqref{BS},\eqref{BS2} the relative integrality conditions are indeed related to the 
the joint asymptotic spectrum of $H$ (the Bohr-Sommerfeld conditions) 
by an embedding $H(X_{\reg})\to \real^n$ which is integral affine 
in leading order.

\end{remark}

\paragraph{Bohr-Sommerfeld conditions}

Similar to the embedding of relative spectra in remark \ref{relAbs}, ``combining" the relative conditions with the formal picture \eqref{diagramm}
through Lemma \eqref{equiIntertwiner} now leads to the following

\begin{suggestion} Let $\star$ be a deformation quantization of a symplectic manifold
$(X,\omega)$ adapted to a Lagrangian submanifold $L$.
Then the formal analogue of (the $\delta$-image of) the pre\-quan\-tum Bohr-\-Som\-mer\-feld class 
\eqref{BS}
is given by
$[\star]_L$.
\end{suggestion}

The point of our  approach to prequantum Bohr-Sommerfeld classes is that 
it reproduces the leading order $\tfrac{1}{\lambda} i^*_L\theta$ of \eqref{BS},\eqref{BS2} already
on the formal algebraic level without involving WKB type methods. This gives
another explanation of the
independence of the leading order conditions on the quantization itself. Moreover,
it generalizes them to arbitrary symplectic manifolds where the integrality condition on relative classes implies one on absolute classes as required by geometric quantization.

On the other hand, the Maslov class is not visible in this approach. 
Let us note though that it can be easily extracted from adapted Fedosov star products as follows:

\paragraph{Maslov index} 
For the Maslov class $\mu$ to be defined consider an additional
  Lagrangian  fibration $\pi$ of $X$ around $L$ (the example in mind is of course the vertical fibration in case of $X=T^*Q$).
  Let $\nabla(\star_F),\nabla(\star'_F)$ be the symplectic connections inside the Fedosov connections of two Fedosov star products $\star_F,\star_F'$ adapted to $L$ and the fibers of  $\pi$ respectively. On $T_LX$, we may further assume these connections to be unitary with respect to some compatible almost $\comp$-structure on $X$, such that the difference
  $\nabla(\star_F) - \nabla(\star_F')$
is identified with a horizontal $\Lie{u}(n)$-valued equivariant 1-form on the unitary frame bundle of $X$ restricted to
 $L$. In this identification
   the Maslov class $\mu$
  defined by $L$ and $\pi$ may be calculated from $\star_F,\star_F'$ as secondary characteristic class:  
\begin{corollary}
$$ \mu = 2i^*_L\tr_{\comp} (\nabla(\star_F) - \nabla(\star_F')).$$
\end{corollary}
  Indeed, by Proposition \ref{adaptFed} the connections are adapted and thus related by local gauge transformations
  $g:  X\supset V \to U(n)$ representing $k$ (see appendix), hence over $L$ we have (cf. \cite{vaisman:1987})
$$
 \tfrac{i}{2\pi}\tr_{\comp}(\nabla-\nabla')= \tfrac{i}{2\pi}tr(g^{-1}dg)
=\tfrac{i}{2\pi}d \ln \det g = 
 \tfrac{1}{2}( \mathrm{det}^2 g)^*d(\ln: \; e^{z}\mapsto z).
$$
 Note that the Liouville class cannot be extracted likewise as characteristic class in general, but if 
 $S =\tau_{\geq 3} H$ and $\nabla_F' = e^{\lad S} \nabla_F e^{-\lad S}$, one calculates
  $ \omega( \nabla -\nabla')  = \delta \tau^3 H  = L_{\ham H}(\nabla)$.
 
\section{Preliminary analogues to symbol calculus of FIOs}

Recall (cf. \cite{bates.weinstein:1995a}) that in generalization of (graphs of) symplectomorphisms a canonical relation $\Lambda$ is defined as Lagrangian submanifold of
$X' \times \bX$, where $\overline{(X,\omega)}:=(X,-\omega)$ denotes the symplectic 
conjugated space. The composition 
$$\Lambda_1 \circ \Lambda_2 
= \Lambda_1 \times_{X'} \Lambda_2
=\pi_{14}(\pi_{12}^{-1} \Lambda_1 \cap \pi_{34}^{-1} \Lambda_2)$$
 (here the $\pi_{ij}$ denote the canonical projections
of $X''\stimes \bX' \stimes X'\stimes \bX$ onto the $i$s and $j$s factor)
may then be identified with the image  of $\Lambda_1\times \Lambda_2$
under symplectic reduction of $X''\stimes \bX' \stimes X'\stimes \bX$ with respect to the
canonical coisotropic manifold $C:=X''\stimes \Delta \stimes \bX$.
$\Lambda_1,\Lambda_2$ are called composable 
if $C$ intersects
$\Lambda_1\times \Lambda_2$ cleanly, then the product is an immersed Lagrangian submanifold $L$. Since multiple points of the immersion correspond to multiple intersections of
$L$ with some fiber of the
characteristic foliation $\pi_{14}|C$ of $C$, it will be an embedding
if the closure of $L$ intersects any fiber at most once.

In terms of the corresponding function algebra $\A:=C^{\infty}(X)$, the above fiber product corresponds either to the 
topological tensor product
$$ \frac{\A'' \hat{\otimes} {\A'}\op}{\I_{\Lambda_1}} \hat{\otimes}_{\A'} \frac{\A' \hat{\otimes} \A\op}{\I_{\Lambda_2}}$$
or, in terms of symplectic reduction, to
$$ N(\I_C)  /  \I_C +\I_{\Lambda_1 \stimes \Lambda_2}, $$
where $N(\I_C)$ is the Poisson normalizer of the vanishing ideal of $C$, which consists
of functions constant along the fibers of the characteristic foliation 
$\pi_{14}| C$.

The deformed analogue of a canonical relation $\Lambda$ will then be
a $\star' \otimes \star\op$-module structure on some flat formal line bundle 
over $\Lambda$, which is of the form $\can \circ S \otimes \phi$
for some $\Lambda$-adapted star product $S(\star' \otimes \star\op)$ and some
flat line bundle $\phi$ over $L$, where the equivalence $S$ is non trivial unless $\Lambda$ is
itself a product. In case of graphs $\Lambda =\Graph \psi$, it was observed by \cite[Prop. 3.1]{bordemann:2004}
that modulo phase these bimodules yield not more than  
homomorphisms of star products deforming $\psi$.

Now, the class of such bimodules is in general unstable under composition,
i.e. modulo phases the tensor products will in general collapse to the classical case.
Indeed, since we cannot expect a complete symbol calculus of quantized symplectomorphisms by 
\cite{vanhove:1951a}, one has to change the category or
select composable objects.

1. The most direct way to do so is to consider composable elements at the
level of 
adapted equivalence classes. 
Given formal deformations $\star_i$ of symplectic manifolds $X_i$, 
composable canonical relations $\Lambda_i\subset X_{i}\stimes \bX_{i+1}$
with quantizations $\ast_i = S_i(\star_{i} \otimes \star_{i+1}\op)$ adapted to 
$\Lambda_i$,
by $i^*_\Delta[\star \otimes \star\op]=0$ one can always find
de Rham representatives of 
$[\ast_1 \otimes \ast_3] $ 
whose restrictions to $C$ are basic, i.e. vanish along the fibers
of $\pi_{14}|C$. 
Suppose this is  true as well for the relative 
class, i.e.
$$ i^*_C([\ast_1 \otimes \ast_2]_{\Lambda_1\stimes \Lambda_3} )\in \pi^*H^2_{dR}(X_1\stimes X_4, \Lambda_1\circ \Lambda_3)\laurent,$$ 
where $\pi:=\pi_{14}|(C,C\cap \Lambda_1\stimes \Lambda_3)$.
Then analogously to \cite[ch. 5]{bordemann:2004} one can 
construct 
an adapted equivalent star product $\star$ which is as well adapted to $C$,
which means that $\bI_C:=\I_C\formal$ is a $\star$-left ideal and its Poisson normalizer $N(\bI_C)$ a $\star$-subalgebra. Then the induced star product 
$(N(\bI_C)/ \bI_C,\ast')$ is adapted to $\Lambda_1 \circ \Lambda_3$
with class determined by
$$i^*_C([\ast_1 \otimes \ast_2]_{\Lambda_1\stimes \Lambda_3} )
= \pi^*([\ast']_{\Lambda_1\circ \Lambda_2} ),$$
since 
$\pi^* \bI_{\Lambda_1 \circ \Lambda_3} = i_C^*(\bI_{\Lambda_1 \stimes \Lambda_3 } \cap N(\bI_C))$.

2. Another strategy is to find some symbol calculus in the derived category of Lagrangian modules. 
The latter was originally considered in the holomorphic case, first for the sheaf of (micro-)differential operators (cf. \cite[ch.5]{kashiwara:2000}) and recently for a 
complex analogue of deformation quantization in \cite{kashiwara.schapira:2005}. 
 However, in the formal real case, 
$\lambda$-convergence problems remain,
this corresponds to the suggestion of Nest and Tsygan in \cite{nest:2004} to modify the
localization procedure.
In case of a cotangent bundle $T^*Q$ one can again
consider fiberwise polynomial algebras over $\comp \laurent$,
then the derived category
of those modules supported on exact sections which 
intersect pairwise transversally 
is stated in \cite{bressler:2002, kontsevich.soibelman:2001} to 
be $A_\infty$ equivalent to the Fukaya category
``quantizing" the Morse complex on $Q$.
This relates $\star$ representations to
mirror symmetry.

\begin{appendix}
\section{Sketch of the symbol calculus of oscillatory distributions}
Nice detailed 
expositions of this theory mainly due to H\"{o}rmander \cite{hoermander:1971} are
\cite{bates.weinstein:1995a},\cite{CDV} for the semiclassical and \cite{duistermaat:1976} for the
conic ($\lambda$-free) case.

 \paragraph{Generating functions} Recall that 
 the image of a section $\eta$ of 
 $T^*Q\stackrel{\pi}{\to}Q$ is Lagrangian if and only if $d\eta=0$, i.e. $\eta$
 is locally the differential of some function on $Q$, since $\eta^*\theta =\eta$ for the canonical form
$\theta= T^*\pi$ on $T^*Q$.  
 A Lagrangian manifold $L\subset T^*Q$ with caustics (i.e. critical values of $\pi | L$)
 now can be locally obtained as well from
 functions $\phi$ on $B:=Q\times \real^k$ (or any surjective submersion $B\stackrel{\rho}{\to} Q$)
 via the image $L_\phi$ of $\im d\phi$ under symplectic reduction of $T^*B$ 
 with respect to the annulator of the vertical bundle $\ker(TB\to \rho^*TQ)^0$,
 which is given in coordinates as
 \begin{equation}\label{eq:LImmersion2}
L_\phi = \left\{ 
\left.\left(q,\piff{ \phi}{ q}(q,\xi)\right)
 \right|
\piff{ \phi}{ \xi}(q,\xi)=0\right \}.
\end{equation}
If $\im d\phi$ and $C$ intersect cleanly,
$L_\phi$ is the immersion $i_\phi$ of the fiber critical set $\Sigma_\phi:=\{\piff{\phi}{\xi}=0\}$ such that caustics correspond to degenerations
$\det \frac{\partial ^2\phi}{\partial \xi^2}=0$, which allows to attack their local classification by considering some equivalent $\phi$ as 
unfolding of $\xi\mapsto \phi(0,\xi)$, cf. \cite{agv:1985}. In fact, as 
proved by H\"{o}rmander, the choice of $\phi$ 
is  locally unique up to (strict) automorphisms of $B$, additions of constants and direct addition of quadratic forms, cf. \cite[Th. 4.18]{bates.weinstein:1995a}. Globally,
\cite{lees:1979} claims that, besides the Liouville class  $[i^*_L\theta]$ occurring already for
Lagrangian sections,
the obstruction to find some $\phi\in C^\infty(Q\times \real^n) $ yielding $L=L_\phi$ is given by 
the $K^1(L)$ homotopy class defined by the ``difference" $k:L\to U(n)/O(n)$ of $TL$ and $T^*_LQ$, as the latter define sections of the bundle of Lagrangian subspaces of $T_LT^*Q$
isomorphic to $L\times U(n)/O(n)$.
In particular, the winding number $\mu.[\gamma]:=\deg [\det^2 \circ k \circ \gamma]$ associated to some 
loop $\gamma: S^1\to L$ is called its Maslov index, and $\mu$ the Maslov class of $L$. 
Thus, for instance, the circle in $\comp$ (harmonic oscillator) doesn't admit a single generating function (note that this obstruction cannot be circumvented by replacing $\real^n$ by
tori, as the first case covers the second in any sense).

Another example to have in mind for the theory of Fourier integral operators is
the classical action functional $S(\gamma)=\int_0^1 L(\dot{\gamma}) dt$
on the fibration which maps a convenient manifold of free paths $[0,1]\to Q$ to their ends in $Q\times Q$. If the corresponding Hamiltonian system given by 
the Legendre transform of $L$ is complete,
 then $L_S$ is (up to symplectic conjugation given by momenta inversion $\overline{(q,p)}=(q,-p)$
 in $T^*Q$) the
graph of its time 1 flow. 
 
 \paragraph{Oscillatory distributions} Basic elements of short wave asymptotics are the
 the WKB-waves, i.e half densities of type 
 $e^{-i S/\lambda}a_\lambda $ with phase $S\in C^\infty(Q)$ and $a_\lambda$ a half density on $Q$ depending polynomially on $\lambda$. Their key property is given by the stationary phase formula: If $\real^k \ni \xi \mapsto \phi(q,\xi)$ is a family of WKB phases such that $L_\phi$ has no caustics, then their superposition
 $  I(\phi,a)(q):=
\int_{\xi\in \real^k}^* e^{i\phi(q,\xi)/\lambda} a(q,\xi) d^k\xi $ is $\lambda$-asymptotically
equivalent to the WKB wave 
\begin{equation}\label{statPhase}
I(\phi,a) \sim  
\left. e^{\frac{i}{\lambda}\phi} e^{-\frac{i\pi}{2}\ind \partial_\xi^2\phi} \left(
\tfrac{a}{ \sqrt{|\det \partial_\xi^2\phi  |}}+ O(\lambda) \right)\right.\circ \rho |_{\Sigma_\phi}^{-1}.
\end{equation}
where $\phi \circ \rho |_{\Sigma_\phi}^{-1}$ coincides with the enveloping phase
(Huygens' principle). It follows that if the composition 
$ \int_Q I(\phi,a) I(\psi,b) $ is well defined,
then its asymptotics may be written as sum (integral) over the intersection points
$L_\phi \cap L_\psi= \{ d(\phi - \psi)=0 \} $ in case of transversal (clean) intersections, 
which allows to lift the singular support of such distributions to their
microsupport $WF(I(\psi,a)):= \supp (a\circ i_\psi^{-1}) \subset \overline{L_\psi}$
in phase space $T^*Q$. 

Indeed, via $\pi|L^{-1}$ the development \eqref{statPhase} may be naturally 
identified with a half density on $L$ such that the singularities of the denominator in
\eqref{statPhase} at caustics appear as artefact of the projection $\pi|L$ onto $Q$. Moreover, if the microsupports of
a set of oscillatory distributions $I(\phi_i,a_i)$ all lie inside some single Lagrangian $L$, then the differences of the pulled back phases in 
 \eqref{statPhase} define locally constant \v{C}ech 1 cocycles on $L$ corresponding to the class \eqref{BS}. In summary, the leading asymptotics of oscillatory distributions $O(L)$ microsupported on $L$ are
described by constant sections of the flat complex line bundle 
$$\left|\wedge\right|^{\half} L \otimes \exp\left( \tfrac{i}{\lambda}i^*_L[\theta] - \tfrac{i\pi}{2}\mu \right)$$ called principal symbols. In particular, the canonical isomorphisms $\overline{O(L)} = O(\overline{L})$ and
$O(L\times L')=O(L)\hat{\otimes} O(L')$ allow us to speak of distributions
microsupported on canonical relations
$L\subset T^*Q\times \overline{T^*Q}$, called  Fourier integral operators (FIOs).
Then the composition of FIOs, if well defined, 
corresponds to 
the  naturally defined composition of principal symbols, cf. \cite[ch.6]{bates.weinstein:1995a}.
  
\paragraph{Reminder on deformation quantization}
In particular, pseudodifferential operators are oscillatory distributions microsupported on 
the identity (i.e.,  the conormal bundle of the diagonal in $T^*(Q\times Q)$), which is naturally identified with
$T^*Q$. Then star products arise as asymptotics of their composition. 
These products were defined purely algebraically
for any symplectic (or Poisson)  manifold $(X,\omega)$ in \cite{bayen.et.al:1978a} (see also the survey articles 
\cite{dito.sternheimer:2002a},\cite{gutt:2000a})
as  a formal deformation 
$(C^\infty(X)\formal ,\star)$
of the classical algebra $(C^\infty(X),\cdot)$ by a sum of bidifferential operators 
\begin{equation}\label{eq:*prod}
\star:= \sum_{i=0}^\infty \lambda^i \star_i 
\end{equation}
with $\star_0=\cdot$, such that 
\begin{enum}
\item $\star$ is associative $[\star,\star]=0$, which may be written order by order as 
\begin{equation}\label{eq:ass}
 2 b \star_n +\sum_{i=1}^{n-1} [\star_i, \star_{n-i}]=0,
\end{equation}
where $b$, $[.,.]$ are the Hochschild coboundary$^\arabic{foot}$
and Gerstenhaber bracket\footnote{
  The Gerstenhaber bracket on $\C$ is the graded supercommutator of the product
   $C\!\circ C'  (f_0,..., f_{k+l}):=
\sum_{i=0}^{k} (-1)^{il}
C(f_0,...,f_{i-1},C'(f_i,...,f_{i+l}),f_{i+l+1},...,f_{k+l})$  turning $\C$ into a differential graded super Lie algebra
w.r.t. the degree $|\C^k|:=k$.
} 
 of $\C(C^\infty(X),C^\infty(X))$, respectively; 
\item the commutator $\frac{i}{\lambda} [.,.]_\star$ deforms the classical Lie structure
    $\frac{i}{\lambda}[f,g]_* =  \{f,g\} \mod O(\lambda)$, thus by \eqref{HKR} we have
\begin{equation}\label{eq:defLie}
\star_1=\tfrac{i}{2} \{.,.\} + b S_1 
\end{equation}
\item $1 \star f= f\star 1= f$.
\end{enum} 
Two such deformations $\star,\star'$ are considered equivalent if they are linked by some algebra isomorphism
$S$ 
 given
by a series of differential operators $S=id + \sum_{i=1}^\infty \lambda^i S_i$; this is
denoted as $\star' = S(\star):= S^{-1}\circ \star \circ S\otimes S $.
In case of a symplectic manifold $(X,\omega)$ the equivalence classes
are in bijection with $\lambda^{-1}[\omega]  + H^2_{dR}(X)\formal$, see 
\cite[App.]{nest.tsygan:1995b} for a short demonstration or \cite{dito.sternheimer:2002a} for more references.

 The equivalence
of two star products $\star,\star'$ in case $H^2_{dR}(X)=0$ was first observed by Lichnerowicz: By applying $1+\lambda S_1$ with $S_1$ as in condition \eqref{eq:defLie}, we may assume
$\star=\star'=\tfrac{i}{2}\{.,.\}$. If now $\star=\star' \mod O(\lambda^k)$, then by 
\eqref{eq:ass} 
\begin{equation}\label{Lichnerowicz1}
b (\star_k-\star_k')=0  \quad \mbox{and} \quad 
b (\star_{k+1}-\star'_{k+1}) + [\star_1, \star_k-\star_k']=0.
\end{equation}
Now, antisymmetrization and $\omega$ provide isomorphisms
\begin{equation}\label{HKR} 
H\C(C^\infty(X); C^\infty(X)) \cong \Gamma(\bigwedge TX) \cong \Omega(X) 
\end{equation}
calculating the Hochschild cohomology$^\arabic{foot}$ of $C^\infty(X)$ (Ko\-stant-\-Hoch\-schild-\-Ro\-sen\-berg (HKR) theorem),
on which $\ad(\star_1)$ acts as de Rham coboundary $d$; hence by $H^2_{dR}(X)=0$
it follows $\star_k-\star_k'=-\ad(\star_1) S_k$ for some derivation
$S_k\in \ker b$. Thus $id + \lambda^k S_k$ provides the induction step 
for constructing an equivalence between $\star$ and $\star'$. 
 
\end{appendix}

\paragraph{Acknowledgements} I would like to thank Stefan Waldmann and Nikolai Neumaier for helpful remarks and in particular, for sharing their expertise in Fedosov's  construction.
Further thanks go to the DFG Graduiertenkolleg ``Physik an Hadronen-Beschleunigern"
for financial support.

\bibliographystyle{halpha}
\bibliography{bibli,dqarticle,dqthesis,dqbook,dqproceeding,dqprocentry}

\end{document}

%% file: commands.tex
\usepackage{amsmath}
\usepackage{amscd}
\usepackage{amssymb}
\usepackage{amsthm}
\usepackage{amsfonts}
\usepackage[latin1]{inputenc}

\theoremstyle{plain}
	\newtheorem{theorem}{Theorem}
	\newtheorem{lemma}{Lemma}
	\newtheorem{corollary}{Corollary}
	\newtheorem{proposition}{Proposition}
	\newtheorem{suggestion}{Suggestion}
\theoremstyle{remark}
	\newtheorem{remark}{Remark}
	\newtheorem*{Remarks}{Remarks}

\newcounter{count}
\newenvironment*{remarks}
{
	\begin{Remarks}
	\begin{list}{\arabic{count}.}
		{\usecounter{count} \setcounter{count}{\value{remark}} 
		\setlength{\leftmargin}{2em} \setlength{\labelsep}{0.5em} 
		\setlength{\itemindent}{0em}\setlength{\labelwidth}{1em}} 
}
{	
	\setcounter{remark}{\value{count}}
	\end{list}
	\end{Remarks}}

\newenvironment{enum}
{
	\begin{list}{\roman{count}.}
	{\usecounter{count} 	
	\setlength{\leftmargin}{2em} \setlength{\labelsep}{0.5em} 
	\setlength{\itemindent}{0em}\setlength{\labelwidth}{1em}		\setlength{\topsep}{0.1ex}\setlength{\parsep}{0ex}\setlength{\itemsep}{0.1ex}} }
{	
	\end{list}}

\DeclareMathOperator{\Aut}{Aut}
\DeclareMathOperator{\Hom}{Hom}

\DeclareMathOperator{\ind}{ind}
\DeclareMathOperator{\tr}{tr}

\DeclareMathOperator{\ad}{ad}

\newcommand{\Pic}{\mathrm{Pic}}
\DeclareMathOperator{\supp}{supp}
\DeclareMathOperator{\Graph}{graph}
\DeclareMathOperator{\im}{im}

\newcommand{\piff}[2]{\frac{\partial#1}{\partial#2} }
\newcommand{\tdiffo}[1]{\left.\tfrac{d}{d#1} \right|_0}

\newcommand{\base}[1]{\frac{\partial}{\partial#1} }

\newcommand{\set}[2]{ \left\{ \left. {#1}  \; \right| \;\: {#2} \right\} }

\newcommand{\ham}{\mathcal{X}}

\newcommand{\Lie}[1]{\mathfrak{#1}}
\newcommand{\real}{\mathbb{R}}
\newcommand{\comp}{\mathbb{C}}

\newcommand{\bij}{\stackrel{\sim}{\to}}
\newcommand{\sur}{\twoheadrightarrow}
\newcommand{\inj}{\hookrightarrow}

\newcommand{\half}{\frac{1}{2}}